\documentclass[11pt, english, draft]{article}
\usepackage{amssymb,amsmath}
\title{\bf An OD-Characterizable Class of Simple Groups}
\author{ {\bf M. Akbari}$^{1}$, {\bf Xiaoyou Chen}$^{2}$ and {\bf Alireza Moghaddamfar}$^{ 3}$\\[0.3cm]
$^1${\em Department of Mathematics, Payame Noor University,}\\ {\em Tehran, Iran,}\\[0.2cm]
 {\em   $^2$College of Sciences, Henan University of Technology,}\\ {\em $450001$, Zhengzhou, China}
\and   and \\[0.2cm]
{\em  $^3$Faculty of Mathematics, K. N. Toosi
University of Technology,}\\
{\em P. O. Box $16315$--$1618$, Tehran, Iran,}\\[0.2cm]
{\em  E-mails}:  {\tt
moghadam@kntu.ac.ir}, { and} {\tt  moghadam@ipm.ir}\\[0.3cm] }

\newtheorem{theorem}{Theorem}[section]
\newtheorem{definition}[theorem]{Definition}

\newtheorem{lm}[theorem]{Lemma}

\begin{document}
\newcommand{\f}{\frac}
\newcommand{\sta}{\stackrel}
\maketitle
\begin{abstract}
\noindent It is proved that  finite nonabelian simple groups $S$  with $\max \pi(S)=37$
are uniquely determined by their order and degree pattern in the class of all finite groups.
 \end{abstract}
\renewcommand{\baselinestretch}{1.1}
\def\thefootnote{ \ }
\footnotetext{{\em AMS subject Classification {\rm 2010}}:
20D05, 20D06, 20D08.\\[0.1cm]
}
\section{Introduction}
 Throughout this note, all the groups under consideration are
{\em finite}, and simple groups are {\em nonabelian}. Given a
group $G$, the {\em spectrum} $\omega(G)$ of $G$ is the set of orders of
elements in $G$. Clearly, the spectrum $\omega(G)$ is {\em closed} and
{\em partially ordered} by the divisibility relation, and hence is
uniquely determined by the set $\mu(G)$ of its elements which are
{\em maximal} under the divisibility relation.

One of the most well-known graphs associated with $G$ is the
{\em prime graph} (or {\em Gruenberg-Kegel graph}) denoted by ${\rm GK}(G)$.
The vertices of  ${\rm GK}(G)$ are the prime divisors of $|G|$ and two distinct
vertices $p$ and $q$ are joined by an edge (written by $p\sim
q$) iff $pq\in\omega(G)$.  If $p_1<p_2<\cdots<p_k$ are all
prime divisors of $|G|$, then we set ${\rm D}(G)=\left(d_G(p_1),
d_G(p_2), \ldots, d_G(p_k )\right)$, where $d_G(p_i)$ denotes the degree
of $p_i$ in the prime graph ${\rm GK}(G)$.  We call this $k$-tuple ${\rm D}(G)$
the {\em degree pattern of} $G$.
 In addition, we denote by
${\rm OD}(G)$ the set of pairwise non-isomorphic finite
groups with the same order and degree pattern as $G$, and put
$h(G)=|{\rm OD}(G)|$. Since there are only
finitely many isomorphism types of groups of order $|G|$,  $1\leqslant h(G)<\infty$.
Now, we have the following definition.
\begin{definition}{\rm
A group $G$ is called {\em $k$-fold OD-characterizable} if $h(G)=k$. Usually, a $1$-fold OD-characterizable group is simply
called {\em OD-characterizable}, and it is called {\em quasi
OD-characterizable} if it is $k$-fold OD-characterizable for some
$k>1$.}
\end{definition}

Notice that OD-characterizability for simple groups $L_2(q)$ was proved in \cite{MZD,  ZS}.
The OD-characterizability problem for alternating groups $A_n$ of degree $n$ $(5\leqslant n\leqslant 100)$
was investigated in \cite{KM}.

Given a prime $p$, $\mathfrak{S}_p$ stands for the set of nonabelian finite simple groups
$S$  such that $p\in \pi (S)\subseteq \{2, 3, 5, \ldots, p\}$.  Based on calculations in the computer algebra system GAP, the sets $\mathfrak{S}_p$ in which $p<10^3$ are determined in \cite{Z}.  According to these results (see also \cite{KM}),
if $S\in \mathfrak{S}_{37}$, then $S$ is isomorphic to one of the following simple groups:
$$L_2(37), \ U_3(11), \ L_2(31^2), \ S_4(31), \ {^2G}_2(27), \ U_3(27), \ L_2(11^3), \ G_2(11), \ U_4(31), \ A_{37}, \ A_{38}, \ A_{39}, \ A_{40}.$$
Previously, it was proved that the following simple groups are OD-characterizable:
$L_2(37)$,  $L_2(31^2)$,  $L_2(11^3)$ \cite{ZS},  $U_3(11)$ \cite{MZD},  ${^2G}_2(27)$ \cite{MZD},
$A_{37}$, $A_{38}$, $A_{39}$, $A_{40}$ \cite{KM}. So, in this note we will concentrate on the OD-characterizability problem for the rest of the groups, and the following is our main result.\\[0.2cm]
{\bf Theorem A.} {\em  The simple groups $S_4(31)$, $U_3(27)$, $G_2(11)$
and $U_4(31)$ are OD-characterizable.}

By combining Theorem A and the above-envisaged results , we obtain the
following corollary.\\[0.2cm]
{\bf Corollary B.} {\em  All simple groups in $\mathfrak{S}_{37}$ are OD-characterizable.}\\[-0.2cm]

We introduce much more notation and definitions (notation used without further explanation is standard).
 Given a group $G$, we denote by $t(G)$ the maximal number of prime divisors of $G$ that are pairwise nonadjacent in  ${\rm GK}(G)$, and by $t(r, G)$ the maximal number of prime divisors of $G$ containing $r$ that are pairwise nonadjacent in  ${\rm GK}(G)$.
Denote by $s(G)$ the number of connected
components of ${\rm GK}(G)$ and by ${\rm GK}_i(G)$, $i=1, 2, \ldots, s(G)$,
the $i$th connected component of ${\rm GK}(G)$. If $G$ is a group of even
order, then we put $2 \in {\rm GK}_1(G)$. It is now easy to see that the order of a group $G$ can be expressed as a product
of some coprime natural numbers $m_i=m_i(G)$, $i=1, 2, \ldots,
s(G)$, with $\pi(m_i)=\pi_i$, where $\pi(m_i)$ signifies the set of all prime divisors of $m_i$.
The numbers $m_1, m_2, \ldots,
m_{s(G)}$ are called the {\it order components of} $G$.

The sequel of this note is organized as follows. In Section 2,
we recall some basic results, especially, on the spectra of
certain finite simple groups, and they will help us find their degree patterns. Section 3 is devoted to the proof
of our main result (Theorem A). Finally, we in Section 4 give a discussion of the relationship
between two groups with the same order and degree pattern.

\section{Preliminaries}
Before proving our main result, we give several lemmas which will be required to determine the degree pattern of the groups under consideration.
\begin{lm}\label{S(4,q)} {\rm (\cite{mazurov-2002})}
Let $q=p^n$, where $p\neq 3$ is a prime. Then, we have
$$\mu(S_4(q))=\left\{(q^2+1)/2, \ (q^2-1)/2, \
p(q+1), \ p(q-1)\right\}.$$
\end{lm}
\begin{lm}\label{spectrumU(3,q)} {\rm (\cite{ZU(3q)})}
If $q$ is a power of an odd prime $p$, then we have:
$$\mu(U_3(q))= \left\{
\begin{array}{ll}
\left\{q^2-q+1, \ q^2-1, \ p(q+1)\right\} & if \ \ \ q \not \equiv -1 \pmod{3},\\[0.3cm]
\left\{(q^2-q+1)/3, \ (q^2-1)/3, \ p(q+1)/3, \
q+1\right\} & if \ \ \ q \equiv -1 \pmod{3}.
\end{array}
\right.$$
\end{lm}
\begin{lm}\label{spectrumG(2,q)} {\rm (\cite{VS})}
If $q$ is a power of a prime $p>5$, then we have:
$$\mu(G_2(q))= \left\{p(q-1), \ p(q+1), \ q^2-1, \ q^2-q+1, \ q^2+q+1
\right\}.$$
\end{lm}
\begin{lm}\label{spectra(u4q)} {\rm (\cite{zav-L4})}
Let $q$ be a power of an odd prime $p$. Denote $d={\rm gcd}(4,
q+1)$. Then $\mu(U_4(q))$ contains the following (and only the
following) numbers:
\begin{itemize}
\item[{\rm (1)}]  $   (q-1)(q^2+1)/d$,
$(q^3+1)/d$, $p(q^2-1)/d$, $q^2-1$;
\item[{\rm (2)}] $p(q+1)$, if and only if $d=4$;
\item[{\rm (3)}] $9$, if and only if $p=3$.
\end{itemize}
\end{lm}

Using Lemmas \ref{spectrumU(3,q)}, \ref{spectra(u4q)}, \ref{S(4,q)} and  \cite[Table 1]{Z},  the required
results concerning some simple groups in $\mathfrak{S}_{37}$ are collected in Table 1.
\begin{center}
\textbf{Table 1.} {\em The orders, spectra and degree patterns of some simple groups in $\mathfrak{S}_{37}$.}\\[0.5cm]
\begin{tabular}{|l|l|l|l|c}
\hline $S$ & $|S|$& $\mu(S)$ & ${\rm D}(S)$ \\[0.1cm]
\hline $S_4(31)$ & $2^{12}\cdot 3^2\cdot 5^2\cdot 13\cdot 31^4\cdot 37$ &
480, 481, 930, 992& (3, 3, 3, 1, 3, 1) \\[0.1cm] $U_3(27)$ &
$2^{5}\cdot 3^9\cdot 7^2\cdot 13\cdot 19\cdot 37$ & 84, 703, 728 & (3, 2, 3, 2, 1,
1) \\[0.1cm]
$G_2(11)$ & $2^{6}\cdot 3^3\cdot 5^2\cdot 7\cdot 11^6\cdot 19\cdot 37$ & 110, 111, 120, 132, 133 & (3,
4, 3, 1, 3, 1, 1)\\[0.1cm]
$U_4(31)$ & $2^{16}\cdot 3^2\cdot 5^2\cdot 7^2\cdot 13\cdot 19\cdot
31^6\cdot 37$ & 992, 960, 7215, 7440, 7448
 & (5, 5, 5, 2, 3, 2, 3, 3) \\[0.1cm] \hline
\end{tabular}
\end{center}

\begin{lm}\label{vasi}{\rm (\cite{vasi})} Let $G$ be a finite group with
$t(G)\geqslant 3$ and $t(2, G)\geqslant 2$, and let $K$ be the
maximal normal solvable subgroup of $G$. Then the quotient group
$G/K$ is an almost simple group, i.e., there exists a non-abelian
simple group $P$ such that $P\leqslant G/K\leqslant {\rm Aut}(P)$.
\end{lm}

\begin{lm}\label{orderout}{\rm (\cite{KM})} Let $S$ be a simple group in  $\bigcup_{5\leqslant p\leqslant 97}\mathfrak{S}_{p}$. Then, we have  $\pi ({\rm Out}(S))\subseteq \{2, 3, 5\}$.
\end{lm}


\section{Proof of the Main Result}
In this section we will prove Theorem A. Before beginning the proof, we draw the prime graphs of the groups $S_4(31)$,
$U_3(27)$, $G_2(11)$ and $U_4(31)$ in Figure  1.


{\small \setlength{\unitlength}{4mm}
\begin{picture}(0,0)(-33,-3)
\put(-29,-5){\circle*{0.35}}
\put(-30,-4){\circle*{0.35}}
\put(-31,-5){\circle*{0.35}}
\put(-30,-6){\circle*{0.35}}
\put(-27,-4){\circle*{0.35}}
\put(-27,-6){\circle*{0.35}}
\put(-31.8,-5.2){\footnotesize 2}%
\put(-30.1,-3.6){\footnotesize 3}%
\put(-28.6,-5.2){\footnotesize 5}%
\put(-30.4,-7){\footnotesize 31}%
\put(-26.65,-4.3){\footnotesize 13}%
\put(-26.65,-6.3){\footnotesize 37}%
\put(-31,-5){\line(1,0){2}}
\put(-31,-5){\line(1,1){1}}
\put(-31,-5){\line(1,-1){1}}
\put(-29,-5){\line(-1,1){1}}
\put(-29,-5){\line(-1,-1){1}}
\put(-27,-6){\line(0,1){2}}
\put(-30,-6){\line(0,1){2}}
\put(-30.5,-8.8){\footnotesize ${\rm GK}(S_4(31))$}


\put(-20,-5){\circle*{0.35}}
\put(-21,-4){\circle*{0.35}}
\put(-22,-5){\circle*{0.35}}
\put(-21,-6){\circle*{0.35}}
\put(-18,-4){\circle*{0.35}}
\put(-18,-6){\circle*{0.35}}
\put(-22.8,-5.2){\footnotesize 2}%
\put(-21.1,-3.6){\footnotesize 3}%
\put(-19.6,-5.2){\footnotesize 7}%
\put(-21.4,-7){\footnotesize 13}%
\put(-17.65,-4.3){\footnotesize 19}%
\put(-17.65,-6.3){\footnotesize 37}%
\put(-22,-5){\line(1,0){2}}
\put(-22,-5){\line(1,1){1}}
\put(-22,-5){\line(1,-1){1}}
\put(-20,-5){\line(-1,1){1}}
\put(-20,-5){\line(-1,-1){1}}
\put(-18,-6){\line(0,1){2}}
\put(-21.5,-8.8){\footnotesize ${\rm GK}(U_3(27))$}


\put(-11,-5){\circle*{0.35}}
\put(-12,-4){\circle*{0.35}}
\put(-13,-5){\circle*{0.35}}
\put(-12,-6){\circle*{0.35}}
\put(-9,-5){\circle*{0.35}}
\put(-7.5,-4){\circle*{0.35}}
\put(-7.5,-6){\circle*{0.35}}
\put(-14.2,-5.2){\footnotesize 11}%
\put(-12.1,-3.6){\footnotesize 2}%
\put(-10.8,-4.6){\footnotesize 3}%
\put(-12.4,-7){\footnotesize 5}%
\put(-9.5,-4.6){\footnotesize 37}%
\put(-7.2,-4){\footnotesize 7}%
\put(-7.2,-6){\footnotesize 19}%
\put(-13,-5){\line(1,0){2}}
\put(-13,-5){\line(1,1){1}}
\put(-13,-5){\line(1,-1){1}}
\put(-11,-5){\line(-1,1){1}}
\put(-11,-5){\line(-1,-1){1}}
\put(-11,-5){\line(1,0){2}}
\put(-12,-6){\line(0,1){2}}
\put(-7.5,-6){\line(0,1){2}}
\put(-12.5,-8.8){\footnotesize ${\rm GK}(G_2(11))$}


\put(-2,-6.25){\circle*{0.35}}%
\put(0,-6.25){\circle*{0.35}}%
\put(2,-6.25){\circle*{0.35}}%
\put(4,-6.25){\circle*{0.35}}%
\put(-2,-4.25){\circle*{0.35}}%
\put(0,-4.25){\circle*{0.35}}%
\put(2,-4.25){\circle*{0.35}}%
\put(4,-4.25){\circle*{0.35}}%
\put(-2.25, -7.25){\footnotesize 7}%
\put(-0.2,-7.25){\footnotesize 2}%
\put(1.8,-7.25){\footnotesize 3}%
\put(3.5,-7.25){\footnotesize 13}%
\put(-2.25, -3.75){\footnotesize 19}%
\put(-0.2,-3.75){\footnotesize 31}%
\put(1.8,-3.75){\footnotesize 5}%
\put(3.5,-3.75){\footnotesize 37}%
\put(-2,-6.25){\line(1,0){6}}
\put(0,-4.25){\line(1,0){4}}
\put(-2,-6.25){\line(0,1){2}}
\put(0,-6.25){\line(0,1){2}}
\put(2,-6.25){\line(0,1){2}}
\put(4,-6.25){\line(0,1){2}}
\put(-2,-4.25){\line(1,-1){2}}
\put(0,-4.25){\line(1,-1){2}}
\put(2,-4.25){\line(1,-1){2}}
\put(0,-6.25){\line(1,1){2}}
\put(2,-6.25){\line(1,1){2}}
\put(-0.7,-8.8){\footnotesize ${\rm GK}(U_4(31))$}
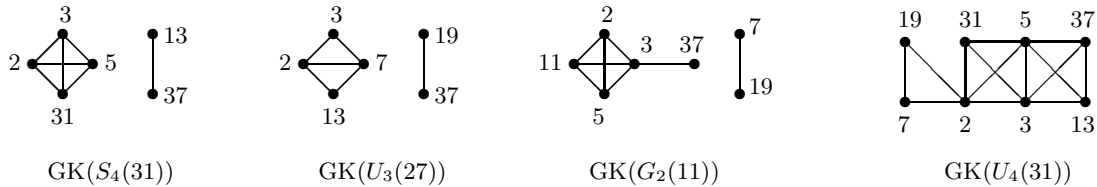
\put(-23,-11){\footnotesize {\bf Fig. 1.} The prime graphs of some
simple groups in $\mathfrak{S}_{37}$.}
\end{picture}}
\vspace{4cm}


\noindent{\em Proof of Theorem A}. Suppose first that $S$ is one of
the simple groups $U_3(27)$, $G_2(11)$
or $U_4(31)$. Let $G$ be a finite group such that $|G|=|S|$
and ${\rm D}(G)={\rm D}(S)$. We have to prove that $G\cong S$. In all these cases
we will prove that $t(G)\geqslant 3$ and $t(2, G)\geqslant 2$.
Therefore, it follows from Lemma \ref{vasi} that there exists a
simple group $P$ such that $P\leqslant G/K\leqslant {\rm
Aut}(P)$, where $K$ is the maximal normal solvable subgroup of
$G$. In addition, we will prove that $P\cong S$, which implies
that $K=1$ and since $|G|=|S|$, $G$ is isomorphic to $S$, as
required. We handle every case singly.

\noindent (a) {\em $S=U_3(27)$}. Let $G$ be a finite group such that
$$|G|=|U_3(27)|=2^{5}\cdot 3^9\cdot 7^2\cdot 13\cdot 19\cdot 37 \ \
\mbox{and} \ \  {\rm D}(G)={\rm D}(U_3(27))=(3, 2, 3, 2, 1, 1).$$
We now consider two cases $19\sim 37$ and $19\nsim 37$, separately.
\begin{itemize}
\item[{\rm (a.1)}] {\em Assume first that $19\sim 37$}. In
this case we immediately have that ${\rm GK}(G)={\rm
GK}(U_3(27))$, and the hypothesis that $|G|=|U_3(27)|$ yields
$G$ and $U_3(27)$ having the same set of order components. Now, by the Main Theorem in
\cite{IKA}, $G$ is isomorphic to $U_3(27)$, as
required.

\item[{\rm (a.2)}] {\em Assume next that $19\nsim 37$}. In this
case, there exists a prime $p\in \pi(G)\setminus \{19, 37\}$ such that $\{p, 19, 37\}$ is an
independent set, otherwise $d_G(19)\geqslant 2$ or $d_G(37)\geqslant 2$, which is impossible. This shows
$t(G)\geqslant 3$.
Moreover, since
$d_G(2)=3$ and $|\pi(G)|=6$, $t(2,G)\geqslant 2$. Thus by Lemma
\ref{vasi} there exists a simple group $P$ such that $P\leqslant
G/K\leqslant {\rm Aut}(P)$, where $K$ is the maximal normal
solvable subgroup of $G$. Let $\pi=\{7, 13, 19, 37\}$. We claim that $K$ is a $\pi'$-group. First of all, if
$\{19, 37\}\subseteq \pi(K)$, then a Hall $\{19, 37\}$-subgroup of
$K$ is an abelian group of order $19\cdot 37$, and hence  $19\sim 37$, which is a contradiction.
Now, assume that  $\{p, q\}=\{
19, 37\}$ and $p$ does not divide the order of $K$ while $q\in \pi(K)$.  Let $Q$ be a Sylow $q$-subgroup of $K$. By Frattini argument $G=KN_G(Q)$. Then,
the normalizer $N_G(Q)$ contains an element of order $p$, say $x$.
Now, $Q\langle x\rangle$ is an abelian group of order $pq$, and so $p\sim q$, again a contradiction.  This shows that $\pi(K)\cap \{19, 37\}=\emptyset$. With the similar arguments, we can verify that  if $13\in \pi(K)$,  then $13$ is adjacent to each of the three vertices 7, 19, and 37, and this forces $d_G(13)\geqslant 3$, which contradicts the hypothesis. Finally, if  $7\in \pi(K)$,  then again $7$ is adjacent to each of the three vertices 13, 19, and 37. Note, however, that  the degree sequence of the subgraph ${\rm GK}(G)\setminus \{7\}$ would be 3, 2, 1,  which is impossible. Therefore, $K$ is a $\pi'$-group.
 Since both $K$ and ${\rm Out}(P)$ are $\pi'$-groups (Lemma \ref{orderout}), $|P|$ is divisible by $7^2\cdot 13\cdot 19\cdot 37$.
Considering the orders of simple groups in
$\mathfrak{S}_{37}$, we conclude that $P$ is isomorphic
to  $U_3(27)$.
Therefore,  $K=1$ and  $G$ is
isomorphic to $U_3(27)$. But then ${\rm GK}(G)={\rm GK}(U_3(27))$ and $19\sim 37$, which is impossible.
\end{itemize}


The proof of the other cases is
quite similar to the proof in the previous case, so we
avoid here full explanation of all details.

\noindent (b) {\em $S=G_2(11)$}. Assume that $G$ is a finite
group such that
$$|G|=|G_2(11)|=2^{6}\cdot 3^3\cdot 5^2\cdot 7\cdot 11^6\cdot 19\cdot 37 \ \ \ {\rm
and} \ \ \ {\rm D}(G)={\rm D}(G_2(11))=(3,
4, 3, 1, 3, 1, 1).$$
We will consider two $7\sim 19$ and $7\nsim 19$, separately.
\begin{itemize}
\item[{\rm (b.1)}] {\em First, suppose that $7\sim 19$}. In
this case,  it follows from ${\rm D}(G)={\rm D}(G_2(11))$ that
the prime graphs of $G$ and $S$ coincide. Thus,
 the hypothesis that $|G|=|G_2(11)|$ yields
$G$ and $G_2(11)$ having the same set of order components. Now, by the Main Theorem in
\cite{ND}, $G$ is isomorphic to $S$, as
required.

\item[{\rm (b.2)}] {\em Next, suppose that $7\nsim 19$}.
It is easy to see that there exists a prime $p\in \pi(G)\setminus \{7, 19\}$
such that  $\{p, 7, 19\}$ is an independent set, and so $t(G)\geqslant 3$.  Moreover,
since $d_G(2)=3$ and $|\pi(G)|=6$, $t(2, G)=3$. Thus by
Lemma \ref{vasi} there exists a simple group $P$ such that
$P\leqslant G/K\leqslant {\rm Aut}(P)$, where $K$ is the maximal
normal solvable subgroup of $G$. Using similar arguments to those
in the previous case, one can show that $K$ is a $\{7,
19, 37\}'$-group and $G$ is isomorphic to $G_2(11)$. But then $7$ is
adjacent to $19$ in ${\rm GK}(G)$, which is a contradiction.
\end{itemize}


\noindent (c) {\em $S=U_4(31)$}. Assume that $G$ is a finite group such
that
$$|G|=|U_4(31)|=
2^{16}\cdot 3^2\cdot 5^2\cdot 7^2\cdot 13\cdot 19\cdot
31^6\cdot 37 \ \ \mbox{and} \ \ {\rm D}(G)={\rm D}(U_4(31))=(5, 5, 5, 2, 3, 2, 3, 3).$$
First of all, we show that $t(G)\geqslant 3$.  To this end, we will consider separately the two cases:
$7\sim 19$ and $7\nsim 19$.  If $7$ is adjacent to $19$ and to another vertex, say $p$, then the induced
graph on $\pi(G)\setminus \{7, 19, p\}$ is not complete, because we have only three vertices with degree $\geqslant 4$.  Therefore, there are at least two nonadjacent vertices $r$ and $s$ in $\pi(G)\setminus \{7, 19, p\}$.  This shows that $\{7, r, s\}$ is an independent set in ${\rm GK}(G)$ and so  $t(G)\geqslant 3$.
If $7$ and $19$ are nonadjacent, then since $d_G(7)=d_G(19)=2$ there exists a vertex which is not adjacent to either of these two vertices, and again we conclude that $t(G)\geqslant 3$.
Moreover, since $d_G(2)=5$ and $|\pi(G)|=8$, $t(2, G)\geqslant 2$.
Thus by Lemma \ref{vasi} there exists a simple group $P$ such
that $P\leqslant G/K\leqslant {\rm Aut}(P)$, where $K$ is the maximal normal solvable subgroup of
$G$.  In addition, $K$ is a $\{7, 19, 37\}'$-group. Indeed, as before, if
$7\in \pi(K)$ or $19\in \pi(K)$, this would yield  $\deg_G(7)\geqslant 3$ or $\deg_G(19)\geqslant 3$,
which is not the case. Finally, if $37\in \pi(K)$, then we obtain $\deg_G(37)\geqslant 4$,  and thus we have a contradiction.
 Since both $K$ and ${\rm Out}(P)$ are $\{7, 19, 37\}'$-groups (Lemma \ref{orderout}), $|P|$ is divisible by $7^2\cdot 19\cdot 37$.
Considering the orders of simple groups in
$\mathfrak{S}_{37}$ yields $P$ isomorphic to
$U_4(31)$. But then $K=1$ and $G$ is isomorphic to $U_4(31)$,
because $|G|=|U_4(31)|$.


Next we concentrate on the simple group $S_4(31)$.

 (d) {\em $S=S_4(31)$}. Suppose that $G$ is a finite group such
that
$$|G|=|S_4(31)|=2^{12}\cdot 3^2\cdot 5^2\cdot 13\cdot 31^4\cdot 37
 \ \ \mbox{and} \ \ {\rm D}(G)={\rm D}(S_4(31))=(3, 3, 3, 1, 3, 1).$$
 We distinguish two cases
separately.
\begin{itemize}
\item[{\rm (d.1)}] {\em Assume first that $13\sim 37$}. In
this case we immediately have that ${\rm GK}(G)={\rm
GK}(S_4(31))$, and since $|G|=|S_4(31)|$ we conclude that
$G$ and $S_4(31)$ have the same set of order components. Now, by the Main Theorem in
\cite{IK}, $G$ is isomorphic to $S_4(31)$, as
required.

\item[{\rm (d.2)}] {\em Assume next that $13\nsim 37$}.  Let $\{p_1, p_2, p_3, p_4\}=\{2, 3, 5, 31\}$.  The prime graph ${\rm GK}(G)$ is depicted in Figure 2. Clearly, $t(G)\geqslant 3$, and since $d_G(2)=3$ and $|\pi(G)|=6$, $t(2, G)\geqslant 2$.

{\small \setlength{\unitlength}{4mm}
\begin{picture}(0,0)(-14,2)
\put(1.8,-1.5){\circle*{0.35}}%
\put(4,0){\circle*{0.35}}%
\put(6,1.3){\circle*{0.35}}%
\put(6,-1.3){\circle*{0.35}}%
\put(8,0){\circle*{0.35}}%
\put(10.3,-1.5){\circle*{0.35}}%
\put(10.5,-2.3){\footnotesize $37$}%
\put(0.7,-2.3){\footnotesize $13$}%
\put(5.75,1.9){\footnotesize $p_4$}%
\put(5.75,-2.1){\footnotesize $p_3$}%
\put(3.3,0.5){\footnotesize $p_1$}%
\put(8.1,0.5){\footnotesize $p_2$}%
\put(1.8,-1.5){\line(3,2){4.2}}
\put(4,0){\line(3,-2){2}}
\put(6,1.3){\line(0,-1){2.5}}
\put(10.3,-1.5){\line(-3,2){4.2}}
\put(8,0){\line(-3,-2){2}}
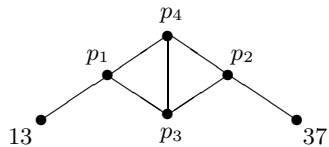
\put(0.5,-4){\footnotesize {\bf Fig. 2.} The prime graph of $G$.}
\end{picture}}
\vspace{3cm}

Thus by Lemma \ref{vasi} there exists a simple group $P$ such
that $P\leqslant G/K\leqslant {\rm Aut}(P)$, where $K$ is the maximal normal solvable subgroup of
$G$.
As before, one can show that $K$ is a
$\{13, 31,  37\}'$-group. Since $K$ and ${\rm Out}(P)$ are $\{13, 31,
37\}'$-groups (Lemma \ref{orderout}), thus $|P|$ is divisible by $13\cdot 31^4\cdot 37$.
Considering the orders of simple groups in
$\mathfrak{S}_{37}$ yields $P$ isomorphic to
$S_4(31)$. But then $K=1$ and $G$ is isomorphic to $S_4(31)$,
because $|G|=|S_4(31)|$.  Therefore ${\rm GK}(G)= {\rm GK}(S_4(31))$
which is disconnected, a contradiction.
\end{itemize}
This completes the proof of theorem. $\Box$

\section{Some Remarks}
Given a finite group $M$, suppose that  $G$ is a finite group with (1) $|G|=|M|$ and (2) ${\rm D}(G)={\rm D}(M)$. In most cases, it follows from the above conditions that they have the same order components. We denote by ${\rm OC}(G)$ the set of order components of $G$. The group $M$ is said to be {\em characterizable by order component}  if, for every finite group $G$, the equality ${\rm OC}(G) = {\rm OC}(M)$ implies the group isomorphism $G\cong M$.  It has already been shown that many simple groups are characterizable by order component (for instance, see \cite{IK, IKA, ND}). Therefore, when under the conditions $|G|=|M|$ and ${\rm D}(G)={\rm D}(M)$ we can conclude that ${\rm OC}(G) = {\rm OC}(M)$, and $M$ is characterizable by order component, it follows that $M$ is  OD-characterizable too. However, in the case when the prime graph of $M$ is connected,  the group $M$ is not necessarily characterizable by order component, but  it may be OD-characterizable. For instance, as we have seen in Theorem A,  the simple group $U_4(31)$
is OD-characterizable, however all nilpotent groups (especially, abelian groups) of order $|U_4(31)|$ have
the same order component, that is $|U_4(31)|$, which means that $U_4(31)$ is not characterizable by order component.

Given a nonnegative integer $n$, we set  $D_n(G)=\{p\in \pi(G) \ | \ d_G(p)=n\}$. Since ${\rm GK}(G)$ is a simple graph,  $D_n(G)=\emptyset$ for all $n\geqslant |\pi(G)|$.  Some information on the prime graph of $G$ is obtained from $D_n(G)$ for some $n$. For instance, since $d_G(p)=0$ if and only if $\{p\}$ is a connected component of ${\rm GK}(G)$, we conclude that $s(G)\geqslant |D_0(G)|$.  On the other hand, if $D_{n-1}(G)\neq \emptyset$, where $n=|\pi(G)|$,  ${\rm GK}(G)$ is connected, that is $s(G)=1$.  In \cite[Theorem B]{suz}), Suzuki proved that if  $L$ is a
finite simple group such that  ${\rm GK}(L)$ is
disconnected, then the connected component ${\rm GK}_i(L)$, $i\geqslant 2$, is a
clique (we recall that a clique is a set of vertices each pair of which is connected by an edge).
As a matter of fact, this is true for
all finite groups not only for finite simple groups. Hence, the prime
graph of an arbitrary finite group $G$ has the following form:
$${\rm GK}(G)=\bigoplus_{i=1}^{s}{\rm GK}_i(G)= {\rm GK}_1(G)\oplus K_{n_2}\oplus \cdots\oplus K_{n_s},$$ where $n_i=|\pi_i(G)|$ $(2\leqslant i\leqslant s)$ and $s=s(G)$. Thus, we conclude that $|D_{n_i-1}(G)|\geqslant n_i$, $(2\leqslant i\leqslant s)$. We denote by
$\pi_i(G)$, $i= 1, 2, \ldots, s(G)$, the set of vertices of
$i$th connected component ${\rm GK}_i(G)$. The
sets $\pi_i(G)$, $i= 1, 2, \ldots, s(G)$, for
finite simple groups $G$ are listed in \cite{k} and \cite{w}.  Under the conditions (1) and (2), if there exists a vertex $p\in D_0(M)$, then $\pi_i(M)=\{p\}=\pi_j(G)$ for some $i, j$.  This restriction helps us determine the group $G$.


\begin{center}
 {\bf Acknowledgments}
\end{center}
This work was done during the second and third authors had a visiting position at the
Department of Mathematical Sciences, Kent State
University, USA. They would like to thank the hospitality of the Department of Mathematical Sciences of KSU.
The second author thanks the funds (2014JCYJ14, 17A110004, 11571129, 11771356).

\end{document}